\documentclass[twocolumn,journal]{IEEEtran}
\usepackage[T1]{fontenc}
\usepackage[latin9]{inputenc}
\usepackage{prettyref}

\usepackage{amsmath, amssymb}
\usepackage{graphicx}
\PassOptionsToPackage{normalem}{ulem}
\usepackage{ulem}
\usepackage{cancel}
\usepackage{tabularx}
\usepackage[unicode=true,
bookmarks=true,bookmarksnumbered=true,bookmarksopen=true,bookmarksopenlevel=1,
 breaklinks=false,pdfborder={0 0 0},pdfborderstyle={},backref=false,colorlinks=false]
 {hyperref}
\hypersetup{pdftitle={Your Title},
 pdfauthor={Your Name},
 pdfpagelayout=OneColumn, pdfnewwindow=true, pdfstartview=XYZ, plainpages=false}

\makeatletter

\usepackage[caption=false,font=footnotesize]{subfig}

\usepackage{cite}
\usepackage{eurosym}

\usepackage{hyperref}

\allowdisplaybreaks

\newtheorem{thm}{Theorem}

\newtheorem{example}{Example}

\def\fnum@figure{Fig.~\thefigure}

\usepackage{multirow,colortbl}
\usepackage{psfrag}

\usepackage{soul}
\soulregister\cite7
\soulregister\ref7
\soulregister\eqref7

\DeclareRobustCommand{\hlyellow}[1]{{\sethlcolor{white}\hl{#1}}} 
\DeclareRobustCommand{\hlcyan}[1]{{\sethlcolor{white}\hl{#1}}} 
\DeclareRobustCommand{\hlred}[1]{{\sethlcolor{white}\hl{#1}}} 
\DeclareRobustCommand{\hlgreen}[1]{{\sethlcolor{white}\hl{#1}}} 
\DeclareRobustCommand{\hlorange}[1]{{\sethlcolor{white}\hl{#1}}} 
\DeclareRobustCommand{\hlpurple}[1]{{\sethlcolor{white}\hl{#1}}} 
\usepackage[colorinlistoftodos]{todonotes}

\newrefformat{fig}{Fig.~\ref{#1}}
\newrefformat{sec}{Section~\ref{#1}}
\newrefformat{subsec}{Section~\ref{#1}}
\newrefformat{app}{Appendix}
\newrefformat{tab}{Table~\ref{#1}}

\newcommand{\etac}{\eta^{\mathrm{C}}}
\newcommand{\etad}{\eta^{\mathrm{D}}}
\newcommand{\cvar}{\overline{c}}
\newcommand{\dvar}{\overline{d}}
\newcommand{\pc}{p^{\mathrm{C}}}
\newcommand{\pd}{p^{\mathrm{D}}}
\newcommand{\rcm}{r^{\mathrm{C}-}}
\newcommand{\rdm}{r^{\mathrm{D}-}}
\newcommand{\rcp}{r^{\mathrm{C}+}}
\newcommand{\rdp}{r^{\mathrm{D}+}}
\newcommand{\pg}{p^{\mathrm{G}}}
\newcommand{\pci}{p^{\mathrm{C}1}}

\newcommand{\pdz}{p^{\mathrm{D}2}}
\newcommand{\PC}{\overline{P^{\mathrm{C}}}}
\newcommand{\PD}{\overline{P^{\mathrm{D}}}}
\newcommand{\PCz}{\overline{P^{\mathrm{C}}_0}}
\newcommand{\PDz}{\overline{P^{\mathrm{D}}_0}}

\newcommand{\eup}{\overline{E}}
\newcommand{\elow}{\underline{E}}
\newcommand{\evar}{\overline{e}}

\newcommand{\redby}{\text{ dominated by }}
\newcommand{\boundon}{\text{ bound on }}
\newcommand{\inch}{\text{ in CH}}

\makeatother

\begin{document}
\title{Tight MIP Formulations for Optimal Operation and Investment of Storage Including Reserves}
\author{$\negthickspace\negthickspace\negthickspace\negthickspace\negthickspace\negthickspace\negthickspace\negthickspace\negthickspace\negthickspace$Maaike B. Elgersma, Germ\'{a}n Morales-Espa\~{n}a, Karen I. Aardal, Niina~Helist\"{o}, Juha~Kiviluoma, and~Mathijs M. de Weerdt
$\negthickspace\negthickspace\negthickspace\negthickspace\negthickspace\negthickspace\negthickspace\negthickspace\negthickspace\negthickspace$}
\maketitle
\begin{abstract}
Fast and accurate large-scale energy system models are needed to investigate the potential of storage to complement the fluctuating energy production of renewable energy systems. However, standard Mixed-Integer Programming (MIP) models that describe optimal investment and operation of these storage units, including the optional capacity to provide up/down reserves, do not scale well. To improve scalability, the integrality constraints are often relaxed, resulting in Linear Programming (LP) relaxations that allow simultaneous charging and discharging, while this is not feasible in practice. To address this, we derive the convex hull of the solutions for the optimal operation of storage for one time period, as well as for problems including investments and reserves, guaranteeing that no tighter MIP formulation or better LP approximation exists for one time period. \hlyellow{When incorporating this convex hull into a multi-period formulation and including it in} large-scale energy system models, the improved LP relaxations can better prevent simultaneous charging and discharging, \hlyellow{and the tighter MIP }\hlpurple{could positively affect}\hlyellow{ the solving time}. We demonstrate this with illustrative case studies of a unit commitment problem and a transmission expansion planning problem. 
\end{abstract}

\begin{IEEEkeywords}
Energy storage systems, mixed-integer linear programming (MIP), linear programming (LP), convex hull, tight formulation, optimal investments, reserves.
\end{IEEEkeywords}

\section*{Nomenclature}
An overview of the notation used throughout this paper is given below.

\vspace{5pt}
\noindent \textit{Parameters:}
\begin{IEEEdescription}[\IEEEusemathlabelsep\IEEEsetlabelwidth{$R^+,\ R^-$}]
\item[$C, D$] Maximum investments of charge/discharge capacities
\item[$E$] Maximum investment of storage capacity
\item[$\elow, \eup$] Minimum/maximum storage capacities
\item[$\eup_0$] Initially installed storage capacity
\item[$\PC, \PD$] Maximum charge/discharge capacities
\item[$\PCz, \PDz$] Initially installed charge/discharge capacities
\item[$R^+, R^-$] Maximum up/down reserve capacities
\item[$\etac, \etad$] Charge/discharge efficiencies 
\item[$\Delta$] Duration of one time period
\item[$\theta$] Ratio between minimum and maximum storage capacity
\end{IEEEdescription}
\noindent \textit{Variables:}
\begin{IEEEdescription}[\IEEEusemathlabelsep\IEEEsetlabelwidth{$R^+,\ R^-$}]
\item[$\cvar, \dvar$] Amount of installed charge and discharge capacity
\item[$e_t$] State of charge at end of time period $t$
\item[$\evar$] Amount of installed storage capacity
\item[$\pc_t, \pd_t$] Amount that is charged and discharged during time period $t$
\item[$r^+_t, r^-_t$] Amount of up and down reserves during time period $t$
\item[$\delta_t$] Binary variable that indicates when the storage unit is charging ($\delta_t=1$) or discharging ($\delta_t=0$)
\end{IEEEdescription}

\section{Introduction}

\subsection{Motivation}

\IEEEPARstart{D}{ue} to climate change, our energy system needs to be decarbonized by 2050. To this end, renewable energy systems need to be implemented on a large scale. However, renewable energy generation fluctuates, due to its intrinsic weather dependency. Storage systems have become a promising solution to complement this fluctuating production, by storing energy when there is an energy surplus, and discharging when there is a deficit. 

To investigate the optimal investment and operation decisions of storage, we need to formulate realistic models for devices capable of storing energy, such as pumped hydro storage, electric vehicles~\cite{srilakshmiEnergyRegulationEV2022}, thermal storage, and some forms of demand response~\cite{singabhattuDistributedEnergyResources2017}. In some cases, we also want to incorporate reserves, which enable more flexibility by representing the possible increase/decrease to the scheduled amount of charged/discharged energy\cite{morales-espanaHiddenPowerSystem2017}. Mixed-Integer Programming (MIP) allows us to correctly model the optimal operation of storage, using both binary and continuous variables. Binary variables are required to correctly model reserves and to ensure that the storage cannot charge and discharge simultaneously~\cite{xuFactoringCycleAging2018,goAssessingEconomicValue2016}.
Constraints in MIP models can describe the characteristics and capabilities of storage units, such as their minimum and maximum state-of-charge, and their (dis)charging capacity. 

\hlyellow{However, there are two big problems that arise when using these storage models to solve large-scale problems. The first is that computational problems arise when MIP storage models are included in large-scale energy system optimization models. MIP models generally do not scale well, since existing algorithms for finding the optimal solution to MIP models run in exponential time (in the worst case). A common practice to reduce the computation time of these large energy models is to }\hlcyan{use a decomposition technique, or to }\hlyellow{relax the integrality restrictions in the MIP formulations. }\hlcyan{In either case, the Linear Programming (LP) relaxation of the storage model is solved, }\hlyellow{which can be solved much faster in practice. However, the problem with these LP }\hlcyan{models for storage operation }\hlyellow{is that their solutions can enforce simultaneous charging and discharging, although this cannot be executed in practice. To address both of these challenges, we improve the LP approximation by tightening the MIP models for optimal storage investment and operation, including reserves.} When incorporated in a large-scale model, these tighter MIP formulations \hlorange{could} speed up the solving time.
When incorporated into a large-scale LP model, the improved LP relaxations can better prevent simultaneous charging and discharging.  

\subsection{Background \& Literature Review}

\hlyellow{Charging and discharging at the same time essentially increases losses, }{\hlcyan{which commonly increases costs\cite{arroyoUseConvexModel2020}. }\hlyellow{Therefore, the optimal solution to the LP relaxation usually does not favor this,}
\hlyellow{but it is beneficial in some cases, such as in the event of negative prices. However, simultaneous charging and discharging can usually not be executed in practice. There are exceptions, such as composite battery systems\cite{elsaadanyLinearModelAggregated2025} or other specific cases with multiple and/or fast storage units. Here, simultaneous charging and discharging might actually be executable by charging a fraction of the time period or storage units, and discharging another fraction\cite{shenModelingArbitrageEnergy2020}. Still, this is not executable when the storage unit is already at maximum capacity, and it is typically not possible in practice at all for generic storage units. Therefore, solutions in which this occurs are unwanted and can lead to large errors.}

\hlyellow{Some examples of situations in which simultaneous charging and discharging in the model solutions can lead to large errors are the following.} Arroyo et al.~\cite{arroyoUseConvexModel2020} present a transmission expansion planning case study, where the optimal investment plan is actually not feasible in practice, because simultaneous charging and discharging is needed to execute it. 
Furthermore, immediate load recovery cannot be guaranteed in demand response models with load shifting when simultaneous up and down shifting of load is allowed~\cite{morales-espanaClassifyingModellingDemand2022}. 
Lastly, significant errors can occur when decomposition algorithms are used to solve large MIP problems faster~\cite{goAssessingEconomicValue2016}. In these decomposition frameworks, the LP relaxation (of the operations/second stage) is solved many times, thus allowing simultaneous charging and discharging and resulting in errors. 

Several attempts to prevent simultaneous charging and discharging in convex problems have been made, such as including pre-contingency operating costs in the objective function~\cite{wenEnhancedSecurityConstrainedUnit2016} or adding conditions for the roundtrip efficiency~\cite{zhaoUsingElectricalEnergy2018}. However, Arroyo et al.~\cite{arroyoUseConvexModel2020} pointed out that it can still occur and illustrated this in two case studies: a unit commitment problem and a transmission expansion planning problem. \hlyellow{Thus, an improved LP-approximation of storage models that can better prevent simultaneous charging and discharging would be a big step toward more useful and better scalable storage models.}

\hlcyan{An improved LP relaxation can also speed up the solving time of MIP storage models. }\hlred{The computational performance of MIP models is highly dependent on their tightness and their size (in terms of the number of constraints and \mbox{variables)\cite{morales-espanaTightCompactMILP2013}}. We say that an MIP model is \textit{tight} if its LP relaxation describes the convex hull of its solutions. We say that MIP model $A$ is \textit{tighter} than another MIP model $B$ if the LP relaxation of MIP model $A$ is contained inside the LP relaxation of MIP model $B$, and it is strictly smaller than the other (so $\text{LP}(A)\subsetneq \text{LP}(B)$). }\hlcyan{Typically, the tightness of MIP models can be improved by adding cuts (additional constraints). However, this also increases the size of the model, which in turn could negatively impact the computational performance.} 

\hlcyan{In related fields, such as unit commitment modeling, much research has been done on finding tighter MIP} \hlorange{\mbox{models\cite{rajanMinimumPolytopesUnit2005,morales-espanaTightCompactMILP2013,gentileTightMIPFormulation2017,huangCuttingPlanesSecurityconstrained2021}}}. These works take into account the important trade-off between model tightness and \hlred{model size.}
\hlyellow{For example, Damci-Kurt et al.\cite{damci-kurtPolyhedralStudyProduction2016} first consider the unit commitment problem with ramping constraints for a single generator in just two time periods. They give a complete linear description of the convex hull of its feasible solutions. This formulation can be incorporated into a multi-generator and multi-period formulation, resulting in a formulation that is tighter than the previously used ramping models, but }\hlred{of similar size. }\hlyellow{They then consider the multi-period case, and find exponential classes of facet-defining inequalities. Though these improve the tightness of the multi-period problem formulation, they }\hlred{heavily increase the size of the model. }\hlyellow{It is only beneficial for the solving time if these constraints are added to the formulation using a separation algorithm, when needed, instead of a priori.}

\hlorange{Tighter formulations for specific storage devices have been found, such as by Baldick et al.\cite{baldickOptimizationFormulationsStorage2023} for pumped storage hydroelectric generators.} For a basic storage operation MIP model, Pozo~\cite{pozoConvexHullFormulations2023} recently obtained two different convex hull formulations: a vertex representation and a hyperplane representation. \hlyellow{Though not clearly stated in their paper, their formulations describe the convex hull of the problem for only one time period. Additional constraints would be needed to describe the convex hull for a multi-period problem, resulting in a model that likely grows exponentially in size w.r.t. the number of time periods. However, the single-period convex hull formulation can be incorporated into a multi-period formulation, which is then tighter than the commonly used formulation.} They demonstrate that incorporating this multi-period LP formulation in a multi-period Set-Point Tracking problem prevents simultaneous charging and discharging more often than the LP relaxation of the commonly used MIP model. \hlorange{Yu et al.\cite{yuTightPowerEnergy2023} presented the same findings, as well as some additional facets for three periods.}

\hlcyan{However, there are several downsides to the LP formulations obtained by Pozo\cite{pozoConvexHullFormulations2023}} \hlorange{and Yu et al.\cite{yuTightPowerEnergy2023}}\hlcyan{. First, we observe that the formulation by Pozo contains several redundant constraints. Half of the constraints can be removed, which }\hlred{decreases the size }\hlcyan{of the model. Second, the convex hull by Pozo formulation }\hlred{does not contain the original binary variable. Therefore, their multi-period LP approximation cannot be used to solve the original multi-period MIP storage operation problem. }\hlorange{Last, neither Pozo nor Yu et al. consider reserves or investment decisions in their models. }

\hlgreen{Thus, developing a tight MIP model for the optimal storage operation problem remains an open question in the literature, as well as finding the convex hull of storage problems that include reserves and investments.}


\subsection{Contributions}

In summary, the main contributions of this paper are the following: \hlcyan{1) We provide a convex hull formulation for optimal storage operation for one period, which consists of fewer constraints than the formulation found by Pozo\cite{pozoConvexHullFormulations2023}. We also provide the corresponding tight MIP model and include it in a multi-period formulation, }\hlred{which can be used to solve the original multi-period storage problem, which is not possible with the LP model by Pozo\cite{pozoConvexHullFormulations2023}. }\hlcyan{2) We provide tighter MIP models for extended problems that include reserves and investment decisions. We prove that their LP-relaxations provide the convex hull of the feasible solutions to the MIP models for one time period, which guarantees that no tighter MIP or better LP approximation for one time period exists. }3) When incorporated in any type of multi-period energy investment/operation system model, we \hlpurple{illustrate} that the LP relaxations of the tighter MIP formulations can decrease the frequency of occurrences of simultaneous charging and discharging.\hlpurple{ } 

The improved MIP formulations and LP relaxations can be used for many different types of planning problems that include storage, as well as for transmission lines. \hlorange{When solving a large-scale MIP problem, tightness of each of the components of the model is an important factor in the solving time. The tighter MIP storage formulations could positively affect the solving time, as shown by Yu et al. \cite{yuTightPowerEnergy2023}.} Furthermore, the improved LP approximation can reduce errors when using decomposition frameworks, as explained earlier.

\section{Background: common storage models}
\label{sec:background}

This section introduces and explains the commonly used MIP models to describe the optimal operation and investment of storage, including reserves. These MIP models are usually incorporated in larger power/energy system models, and the objective depends on the context of these larger problems. Since the objective function is not relevant for obtaining the convex hull of the feasible solutions to this model, it is not included in the model representations below.

\subsection{Modeling Storage Operation Problems}
\label{sec:modelingstorage}

Constraints \eqref{eq:bso-a}-\eqref{eq:bso-f} are commonly used to describe the operation of storage in an MIP model, hereafter named the Basic Operation MIP model (\hyperref[eq:bso]{BO-MIP} model). Here, the index $t$ indicates the time period, and the decision variable $e_t$ represents the state of charge of the storage system \hlyellow{at the end of} time period $t$, which is tracked in \eqref{eq:bso-a}. Decision variables $\pc_t$ and $\pd_t$ represent the amount that is charged and discharged in time period $t$ (with duration $\Delta$), respectively, and parameters $\etac$ and $\etad$ are their corresponding efficiencies. The capacity limits of the storage unit are tracked in \eqref{eq:bso-b}, where parameters $\elow$ and $\eup$ are the min. and max. storage capacities, respectively. The upper and lower bounds for charge and discharge are imposed in \eqref{eq:bso-c} and \eqref{eq:bso-d}. Parameters $\PC$ and $\PD$ are the maximum charge and discharge capacities, respectively.

\noindent\hrulefill \\
\textbf{BO-MIP model:} Basic Operation MIP \\
\vspace{-5pt}
\noindent\hrule
\begin{subequations}
\label{eq:bso}
\begin{align}
&e_{t}  =e_{t-1}+\etac\pc_t\Delta -\frac{1}{\etad}\pd_t\Delta\qquad&\forall t\in T\label{eq:bso-a}\\
&\elow \leq e_{t} \leq\eup \qquad&\forall t\in T\label{eq:bso-b}\\
&\pc_t \leq \PC\delta_t\qquad&\forall t\in T\label{eq:bso-c}\\
&\pd_t \leq \PD(1-\delta_t)\qquad&\forall t\in T\label{eq:bso-d}\\
&e_t,\ \pc_t,\ \pd_t \in \mathbb{R}_{\geq0} \qquad &\forall t\in T\label{eq:bso-e}\\
&\delta_{t}  \in\left\{ 0,1\right\}\qquad&\forall t \in T\label{eq:bso-f}
\end{align}
\end{subequations}
\vspace{-10pt}
\noindent\hrule 
\vspace{10pt} 

Constraints \eqref{eq:bso-c} and \eqref{eq:bso-d} also impose the mutually exclusive condition of charging \textit{or} discharging in each time period, since the binary variable $\delta_{t}$ indicates whether the storage unit is charging ($\delta_{t}=1$) or discharging ($\delta_{t}=0$) in time period $t$. Constraint \eqref{eq:bso-e} defines the variables $e_t$, $\pc_t$ and $\pd_t$ as non-negative continuous variables, and \eqref{eq:bso-f} defines the variable $\delta_{t}$ as binary. Note that to obtain the LP-relaxation of the \hyperref[eq:bso]{BO-MIP} model, hereafter called the \hyperref[eq:bso]{BO-LP} model, the integrality constraint \hlred{on $\delta_t$ is relaxed (\eqref{eq:bso-f} is removed). So $\delta_t$ becomes a continuous variable, which is automatically bounded between 0 and 1 by \eqref{eq:bso-c} and \eqref{eq:bso-d}.}

\hlred{We consider a small numerical example of one time period to illustrate the problem with the \mbox{\hyperref[eq:bso]{BO-LP}} model, which motivates this research. }
\begin{example}
\label{example1}
    \hlyellow{Let us consider a storage unit with $\elow=$ 5MWh, $\eup=$50MWh, $\PC=\PD=$10MW, and $\etac=\etad=$ 0.5 (for simple computational purposes). Now suppose the unit is already at maximum capacity in a certain time period $t^*-1$, so $e_{t^*-1}=50$. In the event of negative prices, or some situation where the surrounding energy system needs to get rid of energy, simultaneous charging and discharging might be beneficial, as explained in the introduction. The \mbox{\hyperref[eq:bso]{BO-LP}} model allows this, so it may output an optimal solutions such as $(e_{t^*-1},\pc_{t^*},\pd_{t^*},\delta_{t^*}) = (50,8,2,0.8)$ for the next time period $t^*$. The state of charge at the end of period $t^*$ is then $e_{t^*}=e_{t^*-1}+\etac\pc_{t^*}\Delta -\frac{1}{\etad}\pd_{t^*}\Delta = 50 + 0.5 \cdot 8 \cdot 1 - 2 \cdot 2\cdot 1 $ $= 50+4-4=50$. So the model allows a solution where it charges as much as it discharges, keeping the battery at the same state of charge, and satisfying all the constraints. However, this can never be executed in practice. Thus, we want to find and LP model that prevents this.}
\end{example}

\subsection{Modeling Storage Operation Problems Including Reserves}
\label{sec:modelingreserves}
As explained earlier, storage can be used to balance the power grid, which is particularly challenging when more renewable energy systems are implemented. Operating reserves, specifically spinning reserves, are useful in this scenario~\cite{morales-espanaHiddenPowerSystem2017}. The up/down reserve provided by a unit is generally defined as the available capacity that the unit can output on top of its scheduled amount in each time period~\cite{morales-espanaMIPFormulationJoint2014}. For storage units specifically, this means that if the storage unit is scheduled to discharge energy, up reserves can be provided by discharging more energy than scheduled, and down reserves by discharging less (or even charging energy). When charging, up reserves can be provided by charging less (or even discharging), and down reserves by charging more.

\noindent\hrulefill \\
\textbf{BOR-MIP model:} Basic Operation incl. Reserves MIP \\
\vspace{-5pt}
\noindent\hrule
\begin{subequations}
\label{eq:bolr}
\begin{align}
&e_{t}  =e_{t-1}+\etac\pc_t\Delta -\frac{1}{\etad}\pd_t\Delta\qquad&\forall t\in T\tag{\ref{eq:bso-a}}\\
&e_{t} \geq \elow + \etac\rcp_t\Delta + \frac{1}{\etad}\rdp_t\Delta \qquad&\forall t\in T\label{eq:bolr-a}\\
&e_{t} \leq\eup -\etac\rcm_t\Delta - \frac{1}{\etad}\rdm_t\Delta \qquad&\forall t\in T\label{eq:bolr-b}\\
&\pc_t + \rcm_t \leq \PC\delta_t\qquad&\forall t\in T\label{eq:bolr-c}\\
&\pd_t + \rdp_t \leq \PD(1-\delta_t)\qquad&\forall t\in T\label{eq:bolr-d}\\
&\pc_t - \rcp_t \geq 0 \qquad&\forall t\in T\label{eq:bolr-e}\\
&\pd_t - \rdm_t \geq 0 \qquad&\forall t\in T\label{eq:bolr-f}\\
&r^-_t = \rcm_t + \rdm_t \leq R^-\qquad&\forall t\in T \label{eq:bolr-g}\\
&r^+_t = \rcp_t + \rdp_t \leq R^+\qquad&\forall t\in T \label{eq:bolr-h}\\
&e_t,\ \pc_t,\ \pd_t,\ \rcp_t,\ \rcm_t,\ \rdp_t,\ \rdm_t \in \mathbb{R}_{\geq0} \quad\hspace{-5pt} &\forall t\in T\label{eq:bolr-i}\\
&\delta_{t}  \in\left\{ 0,1\right\} \qquad &\forall t\in T \tag{\ref{eq:bso-f}}
\end{align}
\end{subequations}
\vspace{-10pt}
\noindent\hrule 
\vspace{10pt}

A common way to model storage operation reserves with an MIP formulation is by adapting the constraints from the \hyperref[eq:bso]{BO-MIP} model, resulting in the Basic Operation incl. Reserves MIP model (\hyperref[eq:bolr]{BOR-MIP} model)\cite{liuSecuredReserveScheduling2021}. Here, \eqref{eq:bolr-a} and \eqref{eq:bolr-b} replace \eqref{eq:bso-b}, since the reserves affect the bounds of the energy storage level and vice versa. Constraint \eqref{eq:bolr-c} replaces \eqref{eq:bso-c}, bounding the amount of energy that can be charged. Constraint \eqref{eq:bolr-c} also bounds the down reserves $\rcm_t$ that can be realized by charging more than planned, and \eqref{eq:bolr-e} bounds the up reserves $\rcp_t$, which can be realized by charging less energy than the planned $\pc_t$. Similarly, \eqref{eq:bolr-d} replaces \eqref{eq:bso-d}, and together with \eqref{eq:bolr-f} imposes bounds on $\pd_t$, $\rdp_t$ and $\rdm_t$. Decision variables $r^-_t$ and $r^+_t$ then represent the total amount of up and down reserves that can be provided in time period $t$. These total reserves are obtained in \eqref{eq:bolr-g} and \eqref{eq:bolr-h} and bounded by some value $R^+$ and $R^-$, respectively.

In this paper, we obtain the convex hull of the \hyperref[eq:bolr]{BOR-MIP} model for one time period. We want to point out that there is an issue with modeling reserves using the \hyperref[eq:bolr]{BOR-MIP} model, namely that it does not fully exploit the flexibility of a storage unit. It only allows limited reserves. We explain this further in Appendix \ref{app:modelingreservesissues}. Here we also point out issues with the formulation presented by Momber et al.\cite{momberPEVStorageMultiBus2014}. This formulation adapts the constraints such that the model does fully exploit the flexibility of the storage unit, but it might not ensure that there is enough capacity to provide the reserves. For those still interested in using the formulation by Momber et al., we provide tight formulations in the online companion\cite[Section 4]{elgersmaOnlineCompanionTight2024}.

\subsection{Modeling Storage Investment Problems Including Reserves}
\label{sec:modelinginvestments}
We can extend the \hyperref[eq:bolr]{BOR-MIP} model, as presented in the previous section, to model investment decisions. The most general way to model the investment of storage is the Basic Investment and operation incl. Reserves MIP model (\hyperref[eq:bgiolr]{BIR-MIP} model)~\cite{morales-espanaImpactLargescaleHydrogen2024}. The operational problem is similar to that in the \hyperref[eq:bolr]{BOR-MIP} model, but it contains additional investment variables $\evar$, $\cvar$, and $\dvar$. The capacity limits of the storage unit, which now depend on the investment variables, are tracked in \eqref{eq:bgiolr-a} and \eqref{eq:bgiolr-b}, where parameter $\eup_0$ is the initially installed maximum storage capacity. Parameter $\theta$ represents the size of the minimum storage capacity as a fraction of the maximum. The upper bounds for charging and discharging and the reserves, also dependent on the investment variables, are imposed in \eqref{eq:bgiolr-c}-\eqref{eq:bgiolr-f}, where parameters $\PCz$ and $\PDz$ are the initially installed charge and discharge capacities, and $C$ and $D$ are the maximum amount of invested charge/discharge capacity. Bounds on the maximum amount of invested capacity are imposed in \eqref{eq:bgiolr-i}-\eqref{eq:bgiolr-k}, where $C$, $D$, and $E$ are typically quite large values.

\noindent\hrulefill \\
\textbf{BIR-MIP model:} Basic Investment and operation incl. Reserves MIP \\
\vspace{-5pt}
\noindent\hrule
\begin{subequations}
\label{eq:bgiolr}
\begin{align}
&e_{t} =e_{t-1}+\etac\pc_t\Delta -\frac{1}{\etad}\pd_t\Delta\qquad&\forall t\in T\tag{\ref{eq:bso-a}}\\
&e_t \geq \theta(\eup_0 + \evar) + \etac\rcp\Delta + \frac{1}{\etad}\rdp\Delta \quad&\forall t\in T \label{eq:bgiolr-a}\\
&e_t \leq \eup_0 + \evar - \etac\rcm\Delta - \frac{1}{\etad}\rdm\Delta \quad&\forall t\in T \label{eq:bgiolr-b}\\
&\pc_t + \rcm_t \leq (\PCz + C)\delta_t  \qquad &\forall t\in T \label{eq:bgiolr-c}\\
&\pc_t +\rcm_t \leq \PCz + \cvar \qquad &\forall t\in T \label{eq:bgiolr-d}\\
&\pd_t + \rdp_t \leq (\PDz + D)(1-\delta_t) \qquad &\forall t\in T \label{eq:bgiolr-e}\\
&\pd_t + \rdp_t \leq \PDz + \dvar \qquad &\forall t\in T \label{eq:bgiolr-f}\\
&\pc_t - \rcp_t \geq 0 \qquad&\forall t\in T\tag{\ref{eq:bolr-e}}\\
&\pd_t - \rdm_t \geq 0 \qquad&\forall t\in T\tag{\ref{eq:bolr-f}}\\
&\cvar \leq C  \label{eq:bgiolr-i}&\\
&\dvar \leq D \label{eq:bgiolr-j}&\\
&\evar \leq E \label{eq:bgiolr-k}&\\
&e_t,\ \pc_t,\ \pd_t,\ \rcp_t,\ \rcm_t,\ \rdp_t,\ \rdm_t \in \mathbb{R}_{\geq0} \quad\hspace{-5pt} &\forall t\in T\tag{\ref{eq:bolr-i}}\\
&\cvar,\ \dvar,\ \evar \in \mathbb{R}_{\geq0} \quad\hspace{-5pt} &\label{eq:bgiolr-m}\\
&\delta_{t}  \in\left\{ 0,1\right\}\qquad &\forall t \in T\tag{\ref{eq:bso-f}}
\end{align}
\end{subequations}
\vspace{-10pt}
\noindent\hrule 
\vspace{10pt}

This model, however, contains many more constraints and variables than the others. A simplified investment and operation model is therefore also commonly used~\cite{tejada-arangoPowerBasedGenerationExpansion2020,goAssessingEconomicValue2016}. We present tight formulations for simplified investment models in the online companion\cite[Section 3]{elgersmaOnlineCompanionTight2024}.

\section{Tight storage models}
\label{sec:results}
In this section, \hlyellow{we explain the methodology and present the resulting tighter MIP formulations for the optimal investment and operation of storage problems (including reserves). In Section \ref{sec:method}, we first explain the outline of the process to obtain and prove the convex hull of the solutions to storage operation problems for one time period. }\hlred{We then apply this to the MIP models introduced in Section \ref{sec:background} to obtain tighter MIP models. }\hlyellow{The full proof for the first presented tight model can be found in Appendix \ref{sec:proofch}. }\hlgreen{The proofs for the other tight models are in the online companion\mbox{\cite[Sections 5 and 6]{elgersmaOnlineCompanionTight2024}}. An overview of all models presented in this paper can also be found in the online companion\mbox{\cite[Section 1]{elgersmaOnlineCompanionTight2024}}, as well as more tight MIP formulations for storage problems with different reserves and investments\mbox{\cite[Sections 2 -- 4]{elgersmaOnlineCompanionTight2024}}.} 


In Section \ref{sec:method}, we also explain what the parameter requirements are for the tighter formulations. Note that if one wants to successfully replace basic storage formulations with the tighter formulations, the parameters might need to be adapted such that they satisfy these requirements. These adaptations are also needed in the case studies, as explained in Appendix \ref{app:cs}.

\subsection{Method: Obtaining the Convex Hull}
\label{sec:method}

The convex hull of an MIP model is a set of constraints that form the tightest possible LP formulation of the solutions to this model. Every vertex of the convex hull is a feasible solution to the MIP model. For more background information on general LP and MIP theory, see~\cite{wolseyIntegerProgramming1998}. In general, we cannot expect to explicitly generate the convex hull of an NP-hard problem, even if we allow for constraint types that contain exponentially many constraints~\cite{papadimitriouComplexityFacetsFacets1982}. However, for storage operation problems in one time period, we are able to find the convex hull by exploiting the disjunctive nature of the MIP. \hlcyan{Balas\cite{balasDisjunctiveProgrammingHierarchy1985} showed how the convex hull of any disjunctive problem can be obtained. This results in a convex hull formulation in a higher dimension than the original disjunctive problem, but it can be projected onto the original dimension.} The obtained projected formulation is still a convex hull\cite[Section 9.2.3]{wolseyIntegerProgramming1998}. Thus, it follows that the obtained formulation describes the convex hull of the original disjunctive problem.

\hlcyan{We can use the theory above to obtain the convex hull of storage problems in the following way.} The first step is to write two disjunctive sets of constraints of the problem for one time period, where one describes the problem when charging, and the other when discharging. For the second step, we can then obtain the convex hull of these two disjunctive sets of constraints, as described Balas~\cite{balasDisjunctiveProgrammingHierarchy1985}. For a detailed explanation of how we did this, we refer to the proof in Appendix \ref{sec:proofch}. However, the obtained convex hull formulation is in a higher dimension than the original formulation. Therefore, in the third step, we project the higher-dimensional formulation to the space of the original variables, by eliminating the extra variables using the Fourier-Motzkin elimination procedure~\cite{fourier1826solution,motzkin1936beitrage}. 

The Fourier-Motzkin elimination procedure eliminates variables by combining all of their upper and lower bounds. The downside of this procedure is that it results in many extra constraints, of which some are facets (constraints describing the convex hull), but many others are redundant. It can be difficult to identify the redundant constraints and prove that they are indeed redundant, since all combinations of constraints need to be checked. \hlred{However, it is typically beneficial to identify and remove redundant constraints, since this could reduce the size of the models. }\hlyellow{For example, for storage operation in one time period, we identified several redundant constraints in the convex hull found by Pozo\cite{pozoConvexHullFormulations2023}. As a result, we found a formulation for the same convex hull consisting of fewer constraints, as proven in Appendix \ref{sec:proofch}.} 

Many constraints are redundant when the parameters satisfy the following reasonable assumptions, which would typically hold in practice. It should always hold that $\PC, R^- \leq \frac{1}{\etac\Delta}(\eup-\elow)$ and $\PD, R^+ \leq \frac{\etad}{\Delta}(\eup-\elow)$, meaning that the charging/discharging capacity and the up/down reserve capacity for one time period cannot be larger than the total storage capacity of the storage unit. Obviously, $\eup > \elow \geq 0$ and $\PC,\PD,R^+,R^-\geq 0$ here. For the efficiencies, $0< \etac, \etad \leq 1$ should hold. Similarly, for the investment problem, it should always hold that $\PCz + C \leq \frac{1}{\etac\Delta}(\eup_0 + E)(1-\phi)$ and $\PDz + D \leq \frac{\etad}{\Delta}(\eup_0 + E)(1-\theta)$, where $0\leq \phi < 1$ and $\PCz, C, \PDz, D, \eup_0, E\geq0$. Thus, to successfully replace some storage formulations with the tighter formulations presented in this paper, the parameters might need to be adapted to satisfy these requirements.

\subsection{Tight MIP Formulation of Storage Operation}
\label{sec:finalmip}
This section presents the Tighter Operation MIP model (\hyperref[eq:chso]{TO-MIP} model) for the optimal operation of storage. Its LP relaxation is obtained by removing the integrality constraint \eqref{eq:bso-f}, meaning $\delta_t$ becomes a continuous variable.

\noindent\hrulefill \\
\textbf{TO-MIP model:} Tighter Operation MIP (tighter version of \hyperref[eq:bso]{BO-MIP} model) \\
\vspace{-5pt}
\noindent\hrule
\begin{subequations}
\label{eq:chso}
\begin{align}
&e_{t}  =e_{t-1}+\etac\pc_t\Delta -\frac{1}{\etad}\pd_t\Delta\qquad&\forall t\in T\tag{\ref{eq:bso-a}}\\
&e_{t-1} \geq \elow + \frac{1}{\etad}\pd_t\Delta  \qquad&\forall t\in T\label{eq:chso-b}\\
&e_{t-1} \leq \eup - \etac\pc_t\Delta  \qquad&\forall t\in T\label{eq:chso-a}\\
&\pc_t \leq \PC\delta_t\qquad&\forall t\in T\tag{\ref{eq:bso-c}}\\
&\pd_t \leq \PD(1-\delta_t)\qquad&\forall t\in T\tag{\ref{eq:bso-d}}\\
&e_t,\ \pc_t,\ \pd_t \in \mathbb{R}_{\geq0} \qquad &\forall t\in T\tag{\ref{eq:bso-e}}\\
&\delta_{t}  \in\left\{ 0,1\right\}\qquad&\forall t \in T\tag{\ref{eq:bso-f}}
\end{align}
\end{subequations}
\vspace{-10pt}
\noindent\hrule 
\vspace{10pt} 

\begin{thm}
    The LP-relaxation of the \hyperref[eq:chso]{TO-MIP} model (\hyperref[eq:chso]{TO-LP}) describes the convex hull of the solutions to the \hyperref[eq:bso]{BO-MIP} model for one time period. \label{thm:chso}
\end{thm}
\begin{IEEEproof}
    See Appendix \ref{sec:proofch}.
\end{IEEEproof}
As a result of Theorem \ref{thm:chso}, the \hyperref[eq:chso]{TO-MIP} model is a tighter MIP model describing the optimal operation of storage than the original \hyperref[eq:bso]{BO-MIP} model. \hlred{We illustrate this in \mbox{Example \ref{example2}}.} The \hyperref[eq:chso]{TO-MIP} model contains the facets \eqref{eq:chso-a} and \eqref{eq:chso-b}, which were also found by Pozo~\cite{pozoConvexHullFormulations2023}. However, the first advantage of the proposed \hyperref[eq:chso]{TO-MIP} model compared to the convex hull (in hyperplane representation) of Pozo is that the proposed model contains the same variables as the original MIP problem, whereas the variable $\delta_t$ is not present in the model by Pozo. Thus, the proposed formulation can \hlred{be used to solve the original multi-period MIP storage operation problem.} The second improvement is that \hlyellow{we paid more attention to removing redundant constraints, which are not facets of the convex hull. As a result}, the \hyperref[eq:chso]{TO-LP} model \hlred{contains fewer constraints} than their convex hull formulation. Note that if we apply Fourier-Motzkin in \hyperref[eq:chso]{TO-MIP} to eliminate variable $\delta_t$, then we can replace the upper bounds on $\pc_t$ and $\pd_t$ in \eqref{eq:bso-c} and \eqref{eq:bso-d} by $\frac{\pc_t}{\PC} + \frac{\pc_t}{\PC} \leq 1\ \forall t\in T$. This results in a model similar to that of Pozo, yet still \hlred{smaller in size} (4 constraints instead of 8). 

\begin{example}
\label{example2}
    \hlyellow{In Section \ref{sec:modelingstorage} we gave a small numerical example (Example \ref{example1}) that illustrates the problem with the \mbox{\hyperref[eq:bso]{BO-LP}} model. The \mbox{\hyperref[eq:bso]{BO-LP}} model allows the solution $(e_{t^*-1},\pc_{t^*},\pd_{t^*},\delta_{t^*}) = (50,8,2,0.8)$, which cannot be executed in practice. However, constraint \eqref{eq:chso-a} in the \mbox{\hyperref[eq:chso]{TO-LP}} model cuts off this solution: $e_{t^*-1}=50$ and $\eup - \etac\pc_{t^*}\Delta = 50 - 0.5 \cdot 8 \cdot 1$ do not satisfy this constraint. This illustrates that the \mbox{\hyperref[eq:chso]{TO-LP}} model is tighter than the \mbox{\hyperref[eq:bso]{BO-LP}} model. Moreover, when considering a problem with more time periods, we can see that the \mbox{\hyperref[eq:chso]{TO-LP}} model is still tighter than the \mbox{\hyperref[eq:bso]{BO-LP}} model. For a two period problem for example, solutions such as $(e_{t^*-1},\pc_{t^*},\pd_{t^*},\delta_{t^*},\ \pc_{t^*+1},\pd_{t^*+1},\delta_{t^*+1}) = $ $(50,8,2,0.8,0,0,0)$, $(50,8,2,0.8,8,2,0.8)$ and $(45,10,0,1,8,2,0.8)$ are all cut off by \mbox{\hyperref[eq:chso]{TO-LP}}, but not by \mbox{\hyperref[eq:bso]{BO-LP}}.}
\end{example}

\subsection{Tight MIP Formulation of Storage Including Reserves}
\label{sec:finalmipreserves}
We now present the Tighter Operation incl. Reserves MIP model (\hyperref[eq:tolr]{TOR-MIP} model) for the optimal operation of storage including limited reserves. 

\noindent\hrulefill \\
\textbf{TOR-MIP model:} Tighter Operation incl. Reserves MIP (tighter version of \hyperref[eq:bolr]{BOR-MIP} model) \\
\vspace{-5pt}
\noindent\hrule
\begin{subequations}
\label{eq:tolr}
\begin{align}
&e_{t}  =e_{t-1}+\etac\pc_t\Delta -\frac{1}{\etad}\pd_t\Delta\qquad&\forall t\in T\tag{\ref{eq:bso-a}}\\
&e_{t-1} \geq \elow + \frac{1}{\etad}\pd_t\Delta + \frac{1}{\etad}\rdp_t\Delta \qquad&\forall t\in T\label{eq:tolr-a}\\
&e_{t-1} \leq\eup -\etac\pc_t\Delta -\etac\rcm_t\Delta \qquad&\forall t\in T\label{eq:tolr-b}\\
&\pc_t + \rcm_t \leq \PC\delta_t\qquad&\forall t\in T\tag{\ref{eq:bolr-c}}\\
&\pd_t + \rdp_t \leq \PD(1-\delta_t)\qquad&\forall t\in T\tag{\ref{eq:bolr-d}}\\
&\pc_t - \rcp_t \geq 0 \qquad&\forall t\in T\tag{\ref{eq:bolr-e}}\\
&\pd_t - \rdm_t \geq 0 \qquad&\forall t\in T\tag{\ref{eq:bolr-f}}\\
&\rcm_t \leq R^-\delta_t\qquad&\forall t\in T \label{eq:tolr-c}\\
&\rdm_t \leq R^-(1-\delta_t)\qquad&\forall t\in T \label{eq:tolr-d}\\
&\rcp_t \leq R^+\delta_t\qquad&\forall t\in T \label{eq:tolr-e}\\
&\rdp_t \leq R^+(1-\delta_t)\qquad&\forall t\in T \label{eq:tolr-f}\\
&r^-_t = \rcm_t + \rdm_t \qquad&\forall t\in T \label{eq:tolr-g}\\
&r^+_t = \rcp_t + \rdp_t \qquad&\forall t\in T \label{eq:tolr-h}\\
&e_t,\ \pc_t,\ \pd_t,\ \rcp_t,\ \rcm_t,\ \rdp_t,\ \rdm_t \in \mathbb{R}_{\geq0} \quad\hspace{-5pt} &\forall t\in T\tag{\ref{eq:bolr-i}}\\
&\delta_{t}  \in\left\{ 0,1\right\} \qquad &\forall t\in T \tag{\ref{eq:bso-f}}
\end{align}
\end{subequations}
\vspace{-10pt}
\noindent\hrule 
\vspace{10pt}

\begin{thm}
    The LP-relaxation of the \hyperref[eq:tolr]{TOR-MIP} model (\hyperref[eq:tolr]{TOR-LP}) describes the convex hull of the solutions to the \hyperref[eq:bolr]{BOR-MIP} model for one time period. \label{thm:chsolr}
\end{thm}
\begin{IEEEproof}
    See the online companion\cite[Section 5]{elgersmaOnlineCompanionTight2024}.
\end{IEEEproof}

As a result of Theorem \ref{thm:chsolr}, the \hyperref[eq:tolr]{TOR-MIP} model is an MIP model describing the optimal operation of storage including reserves, which is tighter than the original \hyperref[eq:bolr]{BOR-MIP} model. Constraints \eqref{eq:bolr-a} and \eqref{eq:bolr-b} are replaced by the facets \eqref{eq:tolr-a} and \eqref{eq:tolr-b}, and \eqref{eq:bolr-g} and \eqref{eq:bolr-h} are replaced by the facets \eqref{eq:tolr-c}-\eqref{eq:tolr-h}.

\subsection{Tight MIP Formulation of Storage Investment Including Reserves}
We now present the tighter models for storage investment and operation, including reserves. The formulations of these tighter models without reserves can simply be obtained by removing the reserve variables. For clarity, these tighter models are given in the online companion\cite[Section 2]{elgersmaOnlineCompanionTight2024}. We start by presenting the more general Tighter Investments and operation incl. Reserves MIP model (\hyperref[eq:tgiolr]{TIR-MIP} model).

\noindent\hrulefill \\
\textbf{TIR-MIP model:} Tighter Investments and operation incl. Reserves MIP (tighter version of \hyperref[eq:bgiolr]{BIR-MIP} model) \\
\vspace{-5pt}
\noindent\hrule
\begin{subequations}
\label{eq:tgiolr}
\begin{align}
&e_{t} =e_{t-1}+\etac\pc_t\Delta -\frac{1}{\etad}\pd_t\Delta\qquad&\forall t\in T\tag{\ref{eq:bso-a}}\\
&e_{t-1} \geq \theta(\eup_0 + \evar) + \frac{1}{\etad}\pd_t\Delta + \frac{1}{\etad}\rdp\Delta \quad&\forall t\in T \label{eq:tgiolr-a}\\
&e_{t-1} \leq \eup_0 + \evar - \etac\pc_t\Delta - \etac\rcm\Delta \quad&\forall t\in T \label{eq:tgiolr-b}\\
&\pc_t + \rcm_t \leq (\PCz + C)\delta_t  \qquad &\forall t\in T \tag{\ref{eq:bgiolr-c}}\\
&\pc_t +\rcm_t \leq \PCz\delta_t + \cvar \qquad &\forall t\in T \label{eq:tgiolr-c}\\
&\pd_t + \rdp_t \leq (\PDz + D)(1-\delta_t) \qquad &\forall t\in T \tag{\ref{eq:bgiolr-e}}\\
&\pd_t + \rdp_t \leq \PDz(1-\delta_t) + \dvar \qquad &\forall t\in T \label{eq:tgiolr-d}\\
&\pc_t - \rcp_t \geq 0 \qquad&\forall t\in T\tag{\ref{eq:bolr-e}}\\
&\pd_t - \rdm_t \geq 0 \qquad&\forall t\in T\tag{\ref{eq:bolr-f}}\\
&\etac\pc_t\Delta \leq (1-\theta)(\eup_0\delta_t + \evar) \qquad &\forall t\in T \label{eq:tgiolr-e}\\
&\frac{1}{\etad}\pd_t\Delta \leq (1-\theta)(\eup_0\delta_t + \evar) \qquad &\forall t\in T \label{eq:tgiolr-f}\\
&\cvar \leq C  \tag{\ref{eq:bgiolr-i}}&\\
&\dvar \leq D \tag{\ref{eq:bgiolr-j}}&\\
&\evar \leq E \tag{\ref{eq:bgiolr-k}}&\\
&e_t,\ \pc_t,\ \pd_t,\ \rcp_t,\ \rcm_t,\ \rdp_t,\ \rdm_t \in \mathbb{R}_{\geq0} \quad\hspace{-5pt} &\forall t\in T\tag{\ref{eq:bolr-i}}\\
&\cvar,\ \dvar,\ \evar \in \mathbb{R}_{\geq0} \quad\hspace{-5pt} &\tag{\ref{eq:bgiolr-m}}\\
&\delta_{t}  \in\left\{ 0,1\right\}\qquad &\forall t \in T\tag{\ref{eq:bso-f}}
\end{align}
\end{subequations}
\vspace{-10pt}
\noindent\hrule 
\vspace{10pt}

\begin{thm}
    The LP-relaxation of the \hyperref[eq:tgiolr]{TIR-MIP} model (\hyperref[eq:tgiolr]{TIR-LP}) describes the convex hull of the solutions to the \hyperref[eq:bgiolr]{BIR-MIP} model for one time period. \label{thm:chgoilr}
\end{thm}
\begin{IEEEproof}
    See the online companion\cite[Section 6]{elgersmaOnlineCompanionTight2024}.
\end{IEEEproof}
As a result of Theorem \ref{thm:chgoilr}, the \hyperref[eq:tgiolr]{TIR-MIP} model is a tighter MIP model describing the optimal investment and operation of storage than the original \hyperref[eq:bgiolr]{BIR-MIP} model. Similar to the case without investments, equations \eqref{eq:bgiolr-a} and \eqref{eq:bgiolr-b} of the \hyperref[eq:bgiolr]{BIR} model are replaced by the facets \eqref{eq:tgiolr-a} and \eqref{eq:tgiolr-b}. Furthermore, equations \eqref{eq:bgiolr-d} and \eqref{eq:bgiolr-f} are replaced by the facets \eqref{eq:tgiolr-c} and \eqref{eq:tgiolr-d} ($\delta_t$ appears in these constraints, making them tighter), and \eqref{eq:tgiolr-e} and \eqref{eq:tgiolr-f} are new facets.

\section{Illustrative case studies}
\label{sec:cs}

To illustrate the improved performance of the newly obtained LP approximations, we apply these formulations to a unit commitment (UC) and a transmission expansion planning (TEP) case study, both presented by Arroyo et al.~\cite{arroyoUseConvexModel2020}. For a detailed description of the case studies, we refer to the paper by Arroyo~\cite{arroyoUseConvexModel2020}. We have slightly adapted these case studies such that they include reserves and investment decisions. \hlcyan{Additionally, we have extended the two-period unit commitment case study to a multi-period problem.} A detailed explanation of the adaptions to case studies is given in Appendix \ref{app:cs}. We have implemented the case studies in Julia using JuMP~\cite{Lubin2023}, and solved them to optimality using the Gurobi solver~\cite[Version 10.0.2]{gurobi} \hlred{with the default settings}. 

\subsection{LP Storage Operation incl. Reserves in Unit Commitment Case Study}
\label{sec:resultscs}

In this section, we illustrate the performance of the proposed multi-period LP relaxation for storage operation (including reserves), and compare it to other LPs by embedding them in a unit commitment (UC) case study. More specifically, we compare the optimal solution of the \textit{basic} MIP storage operation models \hyperref[eq:bso]{BO-MIP} and \hyperref[eq:bolr]{BOR-MIP}, and of their LP relaxations, to the optimal solution of our LP relaxations \hyperref[eq:chso]{TO-LP} and \hyperref[eq:tolr]{TOR-LP}. To this end, we replaced constraints (11)-(15) in the formulation by Wen et al.~\cite{wenEnhancedSecurityConstrainedUnit2016} in the case study by Arroyo et al.~\cite{arroyoUseConvexModel2020}, which describe the basic operation of a storage unit, by the tighter storage operation formulations (with reserves) from this paper.

The unit commitment case study describes an optimal operation problem of two generators and a storage unit over two hourly periods with demands equal to 10MW and 36MW, respectively. The objective is to minimize the operational costs of the generators and the storage unit.

\begin{table}[!h]
    \centering
    \caption{Unit Commitment - Optimal Solutions}
    \resizebox{\columnwidth}{!}{%
    \begin{tabular}{c|cc|cc|cc|cc}
         & \multicolumn{2}{c|}{\hyperref[eq:bso]{BO-MIP}} & \multicolumn{2}{c|}{\hyperref[eq:bso]{BO-LP}} & \multicolumn{2}{c|}{\hyperref[eq:chso]{TO-MIP}} & \multicolumn{2}{c}{\hyperref[eq:chso]{TO-LP}} \\
         \cline{2-9}
         & \multicolumn{2}{c|}{Hour} & \multicolumn{2}{c|}{Hour} & \multicolumn{2}{c|}{Hour} & \multicolumn{2}{c}{Hour} \\
         \cline{2-9}
         Variable & 1 & 2 & 1 & 2 & 1 & 2 & 1 & 2\\
         \hline
         $\pg_{1t}$ (MW) & 12.3 & 27.3 & 13.8 & 28.8 & 12.3 & 27.3 & 12.3 & 27.3\\
         $\pg_{2t}$ (MW) & 0.0 & 2.4 & 0.0 & 0.0 & 0.0 & 2.4 & 0.0 & 2.4\\
         $\pc_t$ (MW) & 2.3 & 0.0 & \textbf{5.8} & 0.0 & 2.3 & 0.0 & 2.3 & 0.0 \\
         $\pd_t$ (MW) & 0.0 & 6.3 & \textbf{2.0} & 7.2 & 0.0 & 6.3 & 0.0 & 6.3 \\
         $e_t$ (MW) & 12.0 & 5.0 & 13.0 & 5.0 & 12.0 & 5.0 & 12.0 & 5.0 \\
         \hline
         Total cost (\$) & \multicolumn{2}{c|}{173.2} & \multicolumn{2}{c|}{130.3} & \multicolumn{2}{c|}{173.2} & \multicolumn{2}{c}{173.2} \\
    \end{tabular}
    }
    \label{tab:uccs}
\end{table}

Table \ref{tab:uccs} shows the optimal solutions to the basic unit commitment case study, as obtained from different formulations for the storage operation part. It can be observed that the optimal solution to the \hyperref[eq:chso]{TO-MIP} model is the same as the optimal solution to the \hyperref[eq:bso]{BO-MIP} model, showing that our MIP formulation represents the exact same problem. More importantly, it can be observed that in the optimal solution to the LP relaxation of this model (the \hyperref[eq:bso]{BO-LP} model), both variables $\pc_t$ and $\pd_t$ have a positive value in the first hour (highlighted in bold), meaning that simultaneous charging and discharging occurs. However, our LP relaxation (\hyperref[eq:chso]{TO-LP}) finds the same solution as the MIP models, so it has successfully prevented simultaneous charging and discharging, as desired.


The storage has a maximum capacity of 13.0MWh. Due to simultaneous charging and discharging in the first hour, energy is lost. This enables generator 1 (with a maximum ramp up rate of 15MW/h) to ramp up its generation ($\pg$) fast enough to satisfy the demand in the second hour. The expensive generator 2 is then not needed, resulting in a cheaper solution. However, simultaneous charging and discharging is not feasible in practice. Thus, the \hyperref[eq:bso]{BO-LP} model finds an operation plan that cannot be implemented.


\begin{table}[!h]
    \centering
    \caption{Unit Commitment Incl. Reserves - Optimal Solutions}
    \resizebox{\columnwidth}{!}{%
    \begin{tabular}{c|cc|cc|cc|cc}
         & \multicolumn{2}{c|}{\hyperref[eq:bolr]{BOR-MIP}} & \multicolumn{2}{c|}{\hyperref[eq:bolr]{BOR-LP}} & \multicolumn{2}{c|}{\hyperref[eq:tolr]{TOR-MIP}} & \multicolumn{2}{c}{\hyperref[eq:tolr]{TOR-LP}} \\
         \cline{2-9}
         & \multicolumn{2}{c|}{Hour} & \multicolumn{2}{c|}{Hour} & \multicolumn{2}{c|}{Hour} & \multicolumn{2}{c}{Hour} \\
         \cline{2-9}
         Variable & 1 & 2 & 1 & 2 & 1 & 2 & 1 & 2\\
         \hline
         $\pg_{1t}$ (MW) & 12.3 & 27.3 & 12.9 & 27.9 & 12.3 & 27.3 & 12.5 & 27.3\\
         $\pg_{2t}$ (MW) & 0.0 & 3.3 & 0.0 & 2.7 & 0.0 & 3.3 & 0.0 & 3.3\\
         $\pc_t$ (MW) & 2.3 & 0.0 & 5.1 & 0.0 & 2.3 & 0.0 & 2.3 & 0.0 \\
         $\pd_t$ (MW) & 0.0 & 5.4 & 2.3 & 5.4 & 0.0 & 5.4 & 0.0 & \textbf{5.4} \\
         $e_t$ (MW) & 12.1 & 6.1 & 12.1 & 6.1 & 12.1 & 6.1 & 12.1 & 6.1 \\
         $\rcp_t$ (MW) & 1.0 & 0.0 & 1.0 & 0.0 & 1.0 & 0.0 & 1.0 & 0.0 \\
         $\rdp_t$ (MW) & 0.0 & 1.0 & 0.0 & 1.0 & 0.0 & 1.0 & 0.0 & 1.0 \\
         $\rcm_t$ (MW) & 1.0 & 0.0 & 1.0 & 1.0 & 1.0 & 0.0 & 1.0 & \textbf{1.0} \\
         $\rdm_t$ (MW) & 0.0 & 1.0 & 0.0 & 0.0 & 0.0 & 1.0 & 0.0 & 0.0 \\
         \hline
         Total cost (\$) & \multicolumn{2}{c|}{191.0} & \multicolumn{2}{c|}{184.1} & \multicolumn{2}{c|}{191.0} & \multicolumn{2}{c}{191.0} \\
    \end{tabular}
    }
    \label{tab:uccsr}
\end{table}

To include reserves, we added constraints that ensure that there are reserves available of at least 1.0MW in each time period. Table \ref{tab:uccsr} shows the different optimal solutions to the unit commitment case study including storage reserves. In the optimal solution to the \hyperref[eq:bolr]{BOR-MIP} model, we can see that the minimum up and down reserves of 1.0MW in both hourly periods are satisfied. In the second hour, more energy is generated by generator 2 and less energy is discharged (compared to the problem without reserves) to achieve this, resulting in a higher total cost. 

It can be observed that our \hyperref[eq:tolr]{TOR-MIP} model results in the same optimal solution as the \hyperref[eq:bolr]{BOR-MIP} model, confirming the correctness of our formulation. In the first hour of the solution to the \hyperref[eq:bolr]{BOR-LP} model, simultaneous charging and discharging occurs again. In the optimal solution to our \hyperref[eq:tolr]{TOR-LP}, simultaneous charging and discharging is prevented. However, it can be observed (highlighted in bold) that down reserves are provided by charging in the second period, though this would not be feasible in practice since the storage is discharging in this period. So the optimal solution to the \hyperref[eq:tolr]{TOR-LP} relaxation still contains some error, but less than for the \hyperref[eq:bolr]{BOR-LP} model. This issue is prevented in the alternative reserves formulation given in Appendix \ref{app:modelingreservesissues} (the \hyperref[eq:bofr]{BOF-MIP} model) since this model does not need a binary variable to model reserves accurately.

\subsection{LP Storage Investment in Transmission Expansion Planning Case Study}

In this section, we illustrate the performance of our multi-period storage investment LP relaxation including reserves, and compare it to other LPs in a transmission expansion planning case study. More specifically, we compare the optimal solution of the \textit{basic} MIP storage investment and operation model 
\hyperref[eq:bgiolr]{BIR-MIP}, and of its LP relaxation, to the optimal solution of our LP relaxation \hyperref[eq:tgiolr]{TIR-LP}. To this end, we replaced constraints (24)-(28) in the formulation presented by Arroyo et al.~\cite{arroyoUseConvexModel2020} (formulation originally by Zhou et al.~\cite{zhaoUsingElectricalEnergy2018}), which describe the basic operation of a storage unit, by the tighter storage investment and operation formulations (with reserves) from this paper. 

The transmission expansion planning problem describes an optimal transmission investment problem of a candidate line of \$5000 between two disconnected buses, with a different generator in each and a storage unit in bus 1. For illustration purposes, the operating condition consists of two hourly periods with both nodal demands respectively equal to 25MW and 100MW. We slightly adapted the case study to also include storage investment decisions and reserves. We substituted the storage capacities by investment variables and changed the objective function to minimize the investment costs of the storage unit, instead of its operational costs. To include reserves, we added constraints that ensure that there are reserves available of at least 10\% of the energy demand in each time period.

Table \ref{tab:uccsi} shows the different optimal solutions to the unit commitment case study with storage investment decisions. Simultaneous charging and discharging occurs in the second hour of the solution to the \hyperref[eq:bgiolr]{BIR-LP} model, but it is prevented by our \hyperref[eq:tgiolr]{TIR-LP} model.

\begin{table}[!h]
    \centering
    \caption{Transmission Expansion Planning - Optimal Solutions}
    \resizebox{\columnwidth}{!}{%
    \begin{tabular}{c|cc|cc|cc|cc}
         & \multicolumn{2}{c|}{\hyperref[eq:bgiolr]{BIR-MIP}} & \multicolumn{2}{c|}{\hyperref[eq:bgiolr]{BIR-LP}} & \multicolumn{2}{c|}{\hyperref[eq:tgiolr]{TIR-MIP}} & \multicolumn{2}{c}{\hyperref[eq:tgiolr]{TIR-LP}} \\
         \cline{2-9}
         & \multicolumn{2}{c|}{Hour} & \multicolumn{2}{c|}{Hour} & \multicolumn{2}{c|}{Hour} & \multicolumn{2}{c}{Hour} \\
         \cline{2-9}
         Variable & 1 & 2 & 1 & 2 & 1 & 2 & 1 & 2  \\
         \hline
         $\pg_{1t}$ (MW) & 135.4 & 145.4 & 125.0 & 115.0 & 135.4 & 145.4 & 135.4 & 145.4\\
         $\pg_{2t}$ (MW) & 0.0 & 0.0 & 25.0 & 100.0 & 0.0 & 0.0 & 0.0 & 0.0 \\
         $\pc_t$ (MW) & 82.9 & 0.0 & 97.5 & \textbf{17.5} & 82.9 & 0.0 & 82.9 & 0.0\\
         $\pd_t$ (MW) & 0.0 & 64.6 & 0.0 & \textbf{12.5} & 0.0 & 64.6 & 0.0 & 64.6\\
         $e_t$ (MW) & 74.6 & 10.0 & 87.8 & 91.1 & 74.6 & 10.0 & 74.6 & 10.0\\
         $\rcp_t$ (MW) & 2.5 & 0.0 & 2.5 & 10.0 & 2.5 & 0.0 & 2.5 & 0.0 \\
         $\rdp_t$ (MW) & 0.0 & 10.0 & 0.0 & 0.0 & 0.0 & 10.0 & 0.0 & 10.0 \\
         $\rcm_t$ (MW) & 2.5 & 0.0 & 2.5 & 10.0 & 2.5 & 0.0 & 2.5 & 0.0 \\
         $\rdm_t$ (MW) & 0.0 & 10.0 & 0.0 & 0.0 & 0.0 & 10.0 & 0.0 & 10.0 \\
         $\cvar$ (MWh) & \multicolumn{2}{c|}{85.4}  & \multicolumn{2}{c|}{100.0}  & \multicolumn{2}{c|}{85.4} & \multicolumn{2}{c}{85.4}  \\
         $\dvar$ (MWh) & \multicolumn{2}{c|}{74.6}  & \multicolumn{2}{c|}{12.5}  & \multicolumn{2}{c|}{74.6} & \multicolumn{2}{c}{74.6}  \\
         $\evar$ (MWh) & \multicolumn{2}{c|}{100.0}  & \multicolumn{2}{c|}{100.0}  & \multicolumn{2}{c|}{100.0} & \multicolumn{2}{c}{100.0}  \\
         \hline
         Total cost (\$) & \multicolumn{2}{c|}{5882} & \multicolumn{2}{c|}{3205} & \multicolumn{2}{c|}{5882} & \multicolumn{2}{c}{5882} \\
    \end{tabular}
    }
    \label{tab:uccsi}
\end{table}

Here, the storage unit in bus 1 has a maximum capacity of 100MWh. Similar to the UC case study, simultaneous charging and discharging in the second hour causes a loss of energy. This allows generator 1 (with a maximum ramp-down rate of 10MW/h) to ramp down fast enough in the second hour, without needing a transmission line to get rid of the excess energy. This results in a cheaper solution, where a line investment is not needed. However, simultaneous charging and discharging is not feasible in practice. Thus, the \hyperref[eq:bgiolr]{BIR-LP} model finds an investment plan that does not allow a feasible operation plan in practice, but our \hyperref[eq:tgiolr]{TIR-LP} model does find a feasible plan.

\subsection{MIP Multi-Period Storage Operation in Unit Commitment Case Study}
\label{sec:resultslargescale}
\hlcyan{In this section, we illustrate the }\hlorange{potential} \hlcyan{performance of the proposed multi-period MIP models for storage operation, and we compare them to basic MIP models by embedding them in a larger multi-period unit commitment case study. More specifically, we compare the computation time that the solver needs to solve the \mbox{\hyperref[eq:bso]{BO-MIP}} model to the time it needs for the \mbox{\hyperref[eq:chso]{TO-MIP}} model. To this end, we extend the original problem with two additional periods, in which the demand is 8.25MW and 0.0MW, respectively. We then repeat this problem 365 times (representing all days in one year), resulting in a problem consisting of 4$\cdot$365 = 1460 time periods. All other parameters remain the same. When solving this multi-period problem using an LP model for storage operation, the optimal solution favors simultaneous charging and discharging in many periods. }\hlgreen{This is for similar reasons as those discussed in Section \ref{sec:resultscs}, namely that charging and discharging simultaneously enables additional ramping, which is more cost-effective than starting up a second, very expensive generator. We now investigate how the models perform.}

\hlcyan{An overview of the results can be found in Table \ref{tab:uclargescale}. }\hlgreen{It can be observed that the \mbox{\hyperref[eq:bso]{BO-MIP}} and \mbox{\hyperref[eq:chso]{TO-MIP}} models again find an optimal solution with the same objective function value, without simultaneous charging and discharging. }\hlcyan{More interestingly, the table displays the total time that is needed to solve the problem to optimality, averaged over 10 runs. Notice that less time is needed to solve the problem using the \mbox{\hyperref[eq:chso]{TO-MIP}} model. More concretely, fewer nodes are explored by the solver in the branch-and-cut algorithm for the \mbox{\hyperref[eq:chso]{TO-MIP}} model. This (partially) explains the difference in solve time. Thus, when solving a large-scale multi-period problem, the tighter \mbox{\hyperref[eq:chso]{TO-MIP}} model }\hlorange{could}\hlcyan{ indeed positively affect the solving time.}

\begin{table}[!h]
    \centering
    \caption{Multi-period Unit Commitment - Solving Details}
    \begin{tabular}{c|c|c|c|c}
         & \hyperref[eq:bso]{BO-MIP} & \hyperref[eq:bso]{BO-LP} & \hyperref[eq:chso]{TO-MIP} & \hyperref[eq:chso]{TO-LP} \\
         \hline
         Total cost (\$) & 70515 & 63053 & 70515 & 63094 \\
         \hline
         Total solve time (s) & 16.0 & <1.0 & \textbf{12.3} & <1.0 \\
         \hline
         \# Explored nodes & 8206 & 1 & \textbf{6950} & 1 \\
         \hline
         \# Periods simultaneous & \multirow{ 2}{*}{0} & \multirow{ 2}{*}{663} & \multirow{ 2}{*}{0} & \multirow{ 2}{*}{\textbf{363}} \\
         charging \& discharging & & & &\\
         \hline
         $\sum_{t\in\mathcal{T}} \pc_t \cdot \pd_t$ & 0 & 3850 & 0 &  \textbf{1403}
    \end{tabular}
    \label{tab:uclargescale}
\end{table}

\hlred{Table \ref{tab:uclargescale} also shows how the \mbox{\hyperref[eq:bso]{BO-LP}} and \mbox{\hyperref[eq:chso]{TO-LP}} models (the LP-relaxations of the MIP models) perform on the multi-period problem. It can be observed that the bound (total cost) provided by the \mbox{\hyperref[eq:chso]{TO-LP}} model is slightly better (closer to the MIP solution), and it prevents simultaneous charging and discharging in many more periods than the \mbox{\hyperref[eq:bso]{BO-LP}} model does. }\hlgreen{The last row gives an indication of how much power is charged and discharged simultaneously in the optimal solution of the models. If this value is zero, no simultaneous charging and discharging occurs. We observe that this value for the \mbox{\hyperref[eq:chso]{TO-LP}} model solution is below 1460, which means that less than 1MW is simultaneously charged and discharged per period on average. This value is more than twice as high for the \mbox{\hyperref[eq:bso]{BO-LP}} model.}

\hlgreen{The formulation by Pozo~\cite{pozoConvexHullFormulations2023} cannot be used to solve this multi-period MIP problem, as their formulation does not contain the original binary variable. Instead, their model would yield results comparable to those of the \mbox{\hyperref[eq:chso]{TO-LP}} model, producing the same solution. However, the multi-period \mbox{\hyperref[eq:chso]{TO-LP}} formulation for this problem includes only half as many constraints as Pozo's multi-period LP formulation. }

\section{Conclusions}
\label{sec:conclusion}
In this paper, we obtained the convex hull of feasible MIP solutions for the optimal operation of storage problem for one time period, as well as for the optimal investment problem, both including reserves. We provide a step by step proof for these convex hull formulations. This guarantees that no tighter MIP formulation exists, meaning no better LP approximation exists for one time period. 

The presented case studies illustrate the improved ability of the LP relaxations of the tighter MIP formulations to prevent simultaneous charging and discharging in a multi-period problem, compared to other LP models. \hlcyan{Furthermore, we illustrated that the tighter MIP formulation }\hlorange{could}\hlcyan{ indeed positively affect the solving time of the model, also for multi-period problems,} \hlorange{though this is case study dependent.}

The improved MIP and LP relaxation can be used as a better proxy of many different types of energy storage systems, as well as transmission lines, in many different large-scale energy system models. Incorporating the improved LP relaxation into large-scale LP models results in more accurate models, that can better prevent simultaneous charging and discharging. Additionally, the tighter MIP formulation has the potential to speed up the solving time of large-scale MIPs. Moreover, it allows the direct application of certain decomposition algorithms, which could further improve the solving time. Thus, the improved models can support the energy transition by accurately modeling the optimal operation of energy storage systems.

\hlyellow{The convex hull for storage optimization models for more than one period remains an open problem. More constraints are needed to describe the convex hull for multiple periods, as we showed in Example \ref{example2}.}\footnote{\hlcyan{We suspect that the convex hull of storage operation grows exponentially in size with respect to the number of time periods. Experiments that we have done using polyhedral geometry software indicate this. Furthermore, by comparing the storage operation problem to the ramping problem for unit commitment, this suspicion is supported.}} Future work could attempt to find facets for these problems in multiple time periods, resulting in a tighter multi-period formulation. The issue with current reserve models that we pointed out in Section \ref{sec:modelingreserves} can also be investigated further. It is desired to have a formulation that accurately models the full flexibility of a storage unit to provide reserves, ideally one that is tighter yet \hlred{small in size}.



\appendix{}

\subsection{Proof of Convex Hull Formulation Optimal Storage Operation in One Time Period}
\label{sec:proofch}
In this section, we obtain the convex hull formulation of the solutions to the basic storage operation problem, as given by the \hyperref[eq:bso]{BO-MIP} model, for one time period. \hlcyan{We use the method as explained in Section \ref{sec:method}, which guarantees that we obtain the convex hull.} We can write the constraints of this model in the following way for one time period $t^*$:
\begin{align*}
&\elow \leq e_{t^*} = e_{t^*-1}+\etac\pc_{t^*}\Delta-\frac{1}{\etad}\pd_{t^*}\Delta \leq\eup \tag{from \ref{eq:bso-a} and \ref{eq:bso-b}}\\
&\elow \leq e_{t^*-1} \leq \eup\tag{from \ref{eq:bso-b}}\\
&\pc_{t^*} \leq \PC\delta_{t^*}\tag{from \ref{eq:bso-c}}\\
&\pd_{t^*} \leq \PD(1-\delta_{t^*})\tag{from \ref{eq:bso-d}}\\
&e_{t^*-1},\ \pc_{t^*},\ \pd_{t^*} \in \mathbb{R}_{\geq0} \tag{from \ref{eq:bso-e}}\\
&\delta_{t^*} \in \{0,1\}. \tag{from \ref{eq:bso-f}}
\end{align*}
We use the method described in \ref{sec:method} to obtain the convex hull. First, we write two disjunctive sets of constraints in Section \ref{sec:disjset1}. We then obtain the convex hull of these two disjunctive sets of constraints in Section \ref{sec:chhighdim1}, as explained by Balas~\cite{balasDisjunctiveProgrammingHierarchy1985}. Lastly, we obtain the convex hull in the dimension of the original formulation by eliminating the extra variable using the Fourier-Motzkin elimination procedure in Section \ref{sec:elime1}.

\subsubsection{Disjunctive constraints sets: charging - discharging}
\label{sec:disjset1}
The two disjunctive sets of constraints of the problem for one time period are given in \eqref{eq:Disjunctive}. We dropped the subscripts $t^*$ and $t^*-1$ from the notation for simplicity. Note that some constraints are redundant. For example, when charging, $e\leq \eup$ is dominated by $e + \etac\pc\Delta \leq \eup$, since the second constraint is a tighter bound on $e$, so it is redundant. 
\begin{align}
\begin{array}{r @{} l}
\text{Charging:} \\
\delta &\ =1\\
\pd &\ =0\\
\cancel{\elow \leq}\ e+\etac\pc\Delta &\ \leq\eup \\
\elow \leq e &\ \cancel{\leq \eup}\\
\pc &\ \leq \PC\\
e &\ \cancel{\geq 0}\\
\pc &\ \geq 0
\end{array} & \quad\mathrm{or}\quad\begin{array}{r@{}l}
\text{Discharging:} \\
\delta &\ =0\\
\pc &\ =0\\
\elow \leq e-\frac{1}{\etad}\pd\Delta &\ \cancel{\leq \eup} \\
\cancel{\elow \leq}\ e &\ \leq \eup\\
\pd &\ \leq \PD\\
e &\ \cancel{\geq 0}\\
\pd &\ \geq 0
\end{array}\label{eq:Disjunctive}
\end{align}

\subsubsection{Convex hull of disjunctive constraint sets}
\label{sec:chhighdim1}
We can then obtain the convex hull of these two disjunctive sets of constraints, as described Balas~\cite{balasDisjunctiveProgrammingHierarchy1985} in the following way. We rename all variables from one set to $\square ^1$ and variables from the other to $\square ^2$, and we multiply all parameters in one set with $\delta^1$, and the others with $\delta^2$. We also need to include the additional constraint $\delta^1 + \delta^2 = 1$. This gives us the constraints in \eqref{eq:chdisjend}, which describe the convex hull of the solutions to the \hyperref[eq:bso]{BO-MIP} model for one time period.
\begin{align}
\begin{split}
e^1 &\geq \elow \delta^1 \\
e^1+\pci\etac & \leq\eup\delta^1 \\
\pci &\leq \PC\delta^1 \\
\pci &\geq 0 \\
e^2-\frac{\pdz}{\etad} &\geq \elow\delta^2 \\
e^2 &\leq \eup \delta^2 \\
\pdz &\leq \PC\delta^2 \\
\pdz &\geq 0 \\
\delta^1 + \delta^2 & = 1 \label{eq:chdisjend}
\end{split}
\end{align}
We can rewrite this problem, such that it looks more similar to the original formulation. We can write $\pci = \pc$ and $\pdz=\pd$. We also rename variable $\delta^1$ as $\delta$, and write $\delta^2$ as $1-\delta$. Furthermore, we let $e=e^1+e^2$, so we can write $e^2$ as $e-e^1$. This results in \eqref{eq:method1}-\eqref{eq:method6}. Constraints \eqref{eq:method3} and \eqref{eq:method6} only contain variables that are also in the \hyperref[eq:bso]{BO-MIP} model, and they are not redundant, thus they are needed in the convex hull (CH) formulation.
\begin{subequations}
\begin{align}
e^1 &\geq \elow \delta \boundon e^1\label{eq:method2}\\
e^1+\etac\pc\Delta & \leq\eup\delta \boundon e^1\label{eq:method1}\\
0 \leq \pc &\leq \PC\delta \inch \label{eq:method3}\\
(e-e^1)-\frac{1}{\etad}\pd\Delta & \geq \elow(1-\delta) \boundon e^1\label{eq:method4}\\
e-e^1 &\leq \eup (1-\delta) \label{eq:method5} \boundon e^1\\
0 \leq \pd &\leq \PD(1-\delta) \inch \label{eq:method6}
\end{align}
\end{subequations}

\subsubsection{Eliminating variable $e^1$}
\label{sec:elime1}
To reduce the dimensionality of the problem, we eliminate variable $e^1$ by applying Fourier-Motzkin elimination. The lower bounds on $e^1$ are \eqref{eq:method2} and \eqref{eq:method5}, which we can rewrite in the following way:
\begin{align*}
e^{1} & \geq \elow \delta \tag{from \ref{eq:method2}}\\
e^{1} & \geq e - \eup (1-\delta). \tag{from \ref{eq:method5}}
\end{align*}
The upper bounds on $e^1$ are \eqref{eq:method1} and \eqref{eq:method4}, which we can rewrite as
\begin{align*}
e^{1} & \leq -\etac\pc\Delta + \eup\delta \tag{from \ref{eq:method1}}\\
e^{1} & \leq e - \frac{1}{\etad}\pd\Delta - \elow (1-\delta). \tag{from \ref{eq:method4}}
\end{align*} 
We now combine the upper and lower bounds on $e^{1}$ according to the Fourier-Motzkin procedure to eliminate the variable. By combining the first lower bound and first upper bound, we obtain
\begin{subequations}
\begin{align}
\elow \delta \leq -\etac\pc\Delta\ + &\ \eup\delta \nonumber\\
\Rightarrow \etac\pc\Delta \leq (\eup - \elow)\delta & \redby \eqref{eq:method3}.\label{eq:red1}
\intertext{Note that this constraint is redundant by \eqref{eq:method3} under reasonable assumptions about the parameter values, as explained in Section \ref{sec:method}. By combining the first lower bound and second upper bound, we obtain}
\elow \delta \leq e\ &- \frac{1}{\etad}\pd\Delta - \elow (1-\delta) \nonumber & \\
\Rightarrow e \geq \elow + \frac{1}{\etad}\pd\Delta & \inch. \label{eq:inchb1}
\intertext{This constraint is not redundant, so it is needed in the convex hull formulation. Following this procedure, we obtain the following constraints by combining the second lower bound and all upper bounds on $e^1$:}
    e \leq \eup - \etac\pc\Delta \label{eq:inchb2} &\inch \\
    \frac{1}{\etad}\pd\Delta \leq (\eup - \elow)(1-\delta) &\redby \eqref{eq:method6}.
\end{align}
\end{subequations}
In conclusion, we have found \eqref{eq:method3}, \eqref{eq:method6}, \eqref{eq:inchb1}, and \eqref{eq:inchb2}, which describe for the convex hull formulation of the solutions to the \hyperref[eq:bso]{BO-MIP} model for one time period. This result is further explained in Section \ref{sec:finalmip}.

\subsection{Case Studies Details}
\label{app:cs}
As explained in Section \ref{sec:cs}, we used the unit commitment and transmission expansion planning case studies by Arroyo et al.~\cite{arroyoUseConvexModel2020}, and adapted them to include reserves and investment decisions. For a full explanation of the case studies, we refer to the original paper. In this section, we explain how we adapted the case studies.

An overview of all the parameter values in the original UC case study, as well as in our adapted case studies (used to obtain the results in Tables \ref{tab:uccs} and \ref{tab:uccsr} can be found in Table \ref{tab:ucparameters}. For the regular optimal operation problem, we slightly adjusted the charging/discharging capacities. As explained at the end of Section \ref{sec:method}, the LP relaxations of our formulations describe the convex hull of the solutions to the MIP models for one time period under the assumption that $\PC \leq \frac{1}{\etac\Delta}(\eup-\elow)$ and $\PD \leq \frac{\etad}{\Delta}(\eup-\elow)$. In the original case study, the values of $\PC$ and $\PD$ do not satisfy this, so we changed these values, as can be seen in Table \ref{tab:csparameters}. Note that this does not affect the optimal mixed-integer solution. 

\begin{table}[h!]
    \centering
    \caption{Storage unit parameter values of case studies}
    \begin{tabular}{c|c|c|c}
        Parameter & Original & Operation & Reserves \\
        \hline
        $\elow$ & 5& 5& 5\\
        $\eup$ & 13& 13& 13\\
        $\PC$ & 12& $\frac{1}{\etac\Delta}(\eup-\elow)$& $\frac{1}{\etac\Delta}(\eup-\elow)$\\
        $\PD$ & 12& $\frac{\etad}{\Delta}(\eup-\elow)$& $\frac{\etad}{\Delta}(\eup-\elow)$\\
        $\etac$ & 0.9& 0.9& 0.9\\
        $\etad$ & 0.9& 0.9& 0.9\\
        $R^+$ & -& -& $\frac{1}{\etac\Delta}(\eup-\elow)$\\
        $R^-$ & -& -& $\frac{\etad}{\Delta}(\eup-\elow)$\\
    \end{tabular}
    \label{tab:ucparameters}
\end{table}

For the case study including reserves, we also needed to satisfy $R^- \leq \frac{1}{\etac\Delta}(\eup-\elow)$ and $R^+ \leq \frac{\etad}{\Delta}(\eup-\elow)$. Thus, we set these parameters equal to these upper limits, as can be seen in Table \ref{tab:csparameters}. Additionally, we added constraints that ensure that the up/down reserves were at least 10\% of the regular energy demand. 

\begin{table}[h!]
    \centering
    \caption{Storage unit parameter values of case studies}
    \begin{tabular}{c|c|c}
        Old parameter & value & New parameter \\
        \hline
        $\elow$ & 5& - \\
        $\eup$ & 100&  $E$ \\
        $\PC$ & 100& $C$ \\
        $\PD$ & 100& $D$ \\
    \end{tabular}
    \label{tab:csparameters}
\end{table}

For the case study including investment decisions, we set the maximum investments of charge/discharge/storage capacity equal to the respective maximum capacities in the other case studies. For $\theta$, we chose to set it to $\elow/\eup$, based on the original case study. The initially installed charge/discharge/storage capacities $\PCz/\PDz/\eup_0$ are set to zero, representing a Greenfield scenario. In the objective function, the investment costs of the storage unit are minimized, instead of the operational costs. Here, the investment costs per MWh charging/discharging capacity are set to \$1/MWh.

\subsection{Modeling Full Flexibility of Reserves Accurately}
\label{app:modelingreservesissues}

We want to point out that there is an issue with modeling reserves using the \hyperref[eq:bolr]{BOR-MIP} model, which is that it does not fully exploit the flexibility of a storage unit. It only allows limited reserves. We illustrate this with an example, sketched in Figure \ref{fig:reservesCD}. Suppose $\PC=\PD=$ 10MW. If the storage is discharging $\pd_t=$ 8MW, then the maximum up reserves the model guarantees is $\rdp_t \leq$ 2MW from \eqref{eq:bolr-d}, which can be realized in real-life by discharging 2MW extra. Also, since the unit is discharging, $\delta_t=$0 forces $\rcp_t, \rcm_t, \pc_t =$0. So the total down reserves in this instance can be $r^-_t=\rdm_t\leq$ 8MW, bounded by \eqref{eq:bolr-f}. However, in real life, we would be able to provide a total of $r^-_t=$ 18MW down reserves, namely by not discharging those 8MW ($\rdm_t=$ 8MW), and by charging 10MW instead ($\rcm_t=$ 10MW). Thus, the \hyperref[eq:bolr]{BOR-MIP} model does not fully exploit the flexibility of the storage unit. In Figure \ref{fig:reservesCD} this missing reserves capacity is indicated by the gray box.

\begin{figure}[h]
    \centering
    \includegraphics[width=0.8\linewidth]{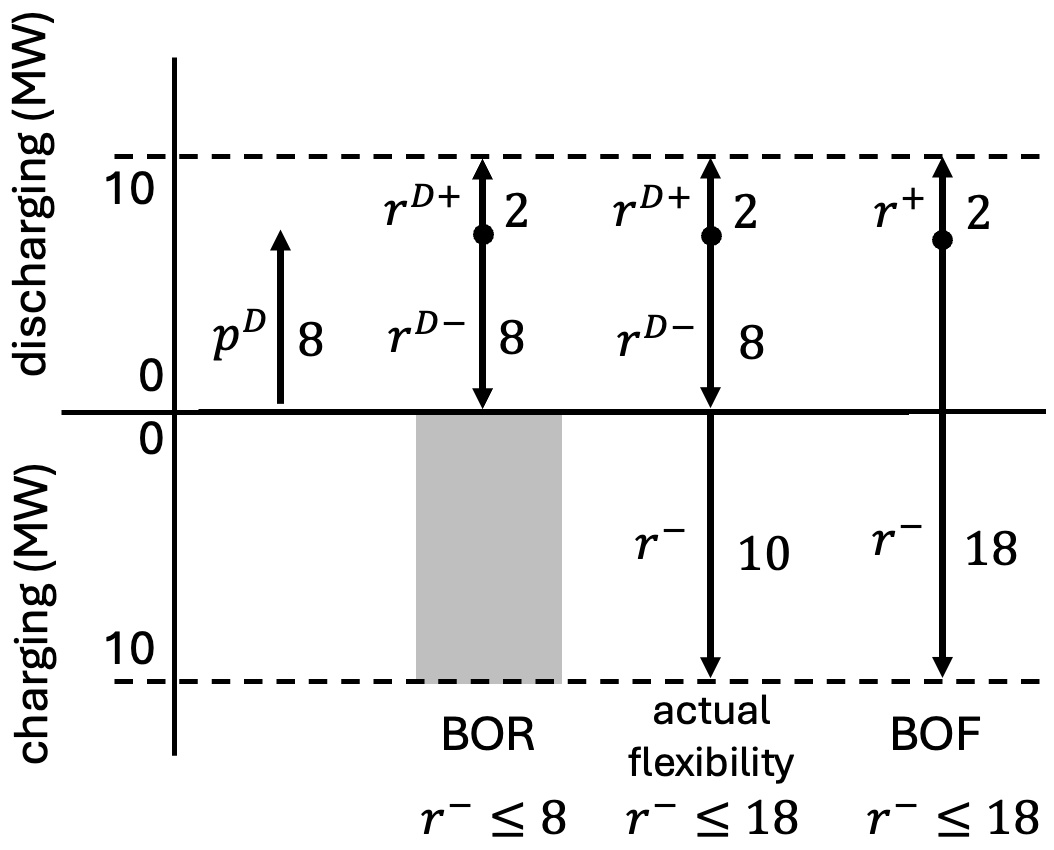}
    \caption{The up and down reserve capacities when discharging, bounded by the charging and discharging capacities, as modeled by different MIP formulations.}
    \label{fig:reservesCD}
\end{figure}

The constraints can be adapted such that the model does fully exploit the flexibility of the storage unit, namely by replacing \eqref{eq:bolr-e} and \eqref{eq:bolr-f} by
\begin{subequations}
    \begin{align}
    \pc_t - \rcp_t \geq -\PD \label{eq:reserves1}\\
    \pd_t - \rdm_t \geq -\PC. \label{eq:reserves2}
\end{align}
\end{subequations}
Constraints \eqref{eq:reserves1} and \eqref{eq:reserves2}  together with \eqref{eq:bolr-c} and \eqref{eq:bolr-d} can be rewritten as \eqref{eq:bofr-b} and \eqref{eq:bofr-c}. This results in the \hyperref[eq:bofr]{BOF-MIP} model (Basic Operation with Flexible Reserves MIP model), as presented by Momber et al.\cite{momberPEVStorageMultiBus2014}. Now, \eqref{eq:bofr-b} and \eqref{eq:bofr-c} ensure that the down reserves $r^-_t$ and up reserves $r^+_t$, respectively, are essentially bounded by $\PD+\PC$. Thus, this formulation does fully exploit the flexibility of the storage unit, as illustrated in Figure \ref{fig:reservesCD}. Lastly, the capacity of the reserves is also bounded by the current state of charge, as modeled by \eqref{eq:bofr-a}. We provide a tighter formulation of this model in the online companion\cite[Section 4]{elgersmaOnlineCompanionTight2024}.

\newpage
\noindent\hrulefill \\
\textbf{BOF-MIP model:} Basic storage operation incl. flexible reserves MIP \\
\vspace{-5pt}
\noindent\hrule
\begin{subequations}
\label{eq:bofr}
\begin{align}
&e_{t}  =e_{t-1}+\etac\pc_t\Delta -\frac{1}{\etad}\pd_t\Delta\qquad&\forall t\in T\tag{\ref{eq:bso-a}}\\
&\elow + \frac{1}{\etad}r^+_t\Delta \leq e_{t}  \leq\eup -\etac r^-_t\Delta \qquad&\forall t\in T\label{eq:bofr-a}\\
&\pc_t \leq \PC\delta_t\qquad&\forall t\in T\tag{\ref{eq:bso-c}}\\
&\pd_t \leq \PD(1-\delta_t)\qquad&\forall t\in T\tag{\ref{eq:bso-d}}\\
&\pc_t-\pd_t+r_{t}^{-}  \leq \PC\qquad&\forall t\in T\label{eq:bofr-b}\\
&-\pc_t+\pd_t+r_{t}^{+}  \leq \PD\qquad&\forall t\in T\label{eq:bofr-c}\\
&r^+_t \leq R^+\qquad&\forall t\in T \label{eq:bofr-d}\\
&r^-_t \leq R^-\qquad&\forall t\in T \label{eq:bofr-e}\\
&e_t,\ \pc_t,\ \pd_t,\ r^+_t,\ r^-_t \in \mathbb{R}_{\geq0} \qquad &\forall t\in T\label{eq:bofr-f}\\
&\delta_{t}  \in\left\{ 0,1\right\} \qquad &\forall t\in T \tag{\ref{eq:bso-f}}
\end{align}
\end{subequations}
\vspace{-10pt}
\noindent\hrule 
\vspace{10pt}

However, we want to point out that \eqref{eq:bofr-a} might not ensure that the reserves are bounded correctly. Thus, it might promise reserves that cannot be realized in practice. Let us consider another example, sketched in Figure \ref{fig:reservesE}. The maximum capacity of the storage unit is 10MWh, the maximum charging/discharging capacities are again 10MW, the charging and discharging efficiencies are 50\% and 100\%, respectively, and $\Delta=$ 1. Suppose the state-of-charge is initially at maximum capacity and the storage is scheduled to discharge 8MW, so the storage level becomes 2MWh. Then \eqref{eq:bofr-b} would allow a down reserve of 18MW, just like in the previous example in Figure \ref{fig:reservesCD}. Due to the total battery capacity and the current storage level, \eqref{eq:bofr-a} would constrain this down reserve such that $\etac r^- \Delta \leq$ 8MW, so $r^-\leq$ 16MW. However, in reality, the down reserve would be provided by not realizing the scheduled discharge, resulting in a down reserve of capacity of 8MW. An error is introduced due to the difference in charging and discharging efficiencies. The promised down reserve of 16MW cannot be provided, showing that \eqref{eq:bofr-a} is not restricted enough to model the problem accurately. 

\begin{figure}[h]
    \centering
    \includegraphics[width=0.8\linewidth]{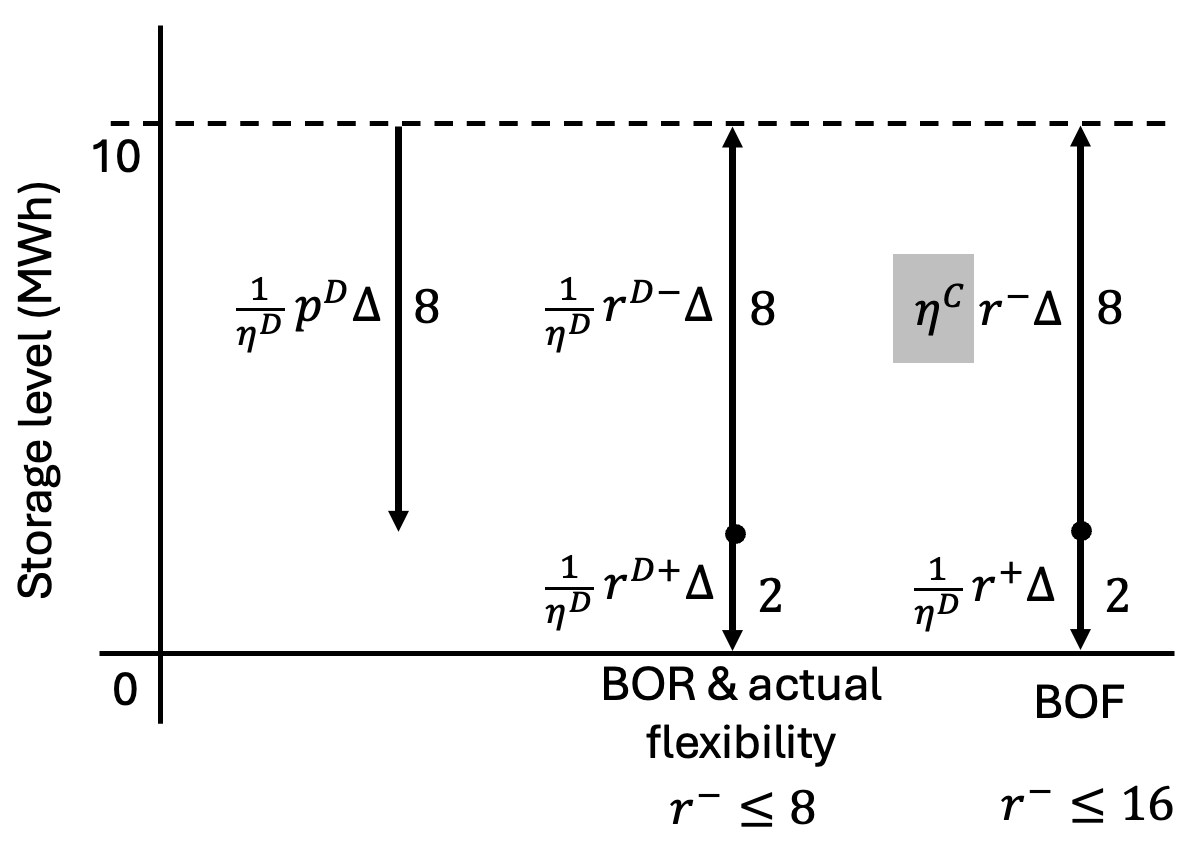}
    \caption{The up and down reserve capacities when discharging, bounded by the total storage capacity, as modeled by different MIP formulations.}
    \label{fig:reservesE}
\end{figure}

\section*{Acknowledgment} 
This research received funding from the Dutch Research Council (NWO), as part of the Energy Systems Integration (ESI-far) program, project number ESI.2019.008, and was supported by eScienceCenter under project number NLeSC C 21.0226. J.~Kiviluoma acknowledges funding from the European Climate, Infrastructure and Environment Executive Agency under the European Union's HORIZON Research and Innovation Actions under grant agreement N101095998. N. Helist\"{o} funding from the European Union's Horizon 2020 research and innovation program under grant agreement 864276.

{\scriptsize{}\bibliographystyle{IEEEtran}
\bibliography{storagelib}
}{\scriptsize\par}

\newpage

\end{document}